\documentclass[preprint,12pt]{elsarticle}

\usepackage{amsmath}
\usepackage{amssymb}
\usepackage{amsthm}

\journal{Games and Economic Behavior}

\newcommand{\hatk}{\hat{k}}
\newcommand{\hatm}{\hat{m}_n}
\newcommand{\hatq}{\hat{q}}
\newcommand{\hatdelta}{\hat{\delta}_n}

\begin{document}

\begin{frontmatter}



\title{The time scales of the aggregate learning and sorting in market entry games with large number of players}


\author{Misha Perepelitsa}

\ead{misha@math.uh.edu}

\address{Department of Mathematics\\
University of Houston\\
4800 Calhoun Rd, Houston, TX, 77096}

\begin{abstract}
We consider the dynamics of player's strategies in repeated market games, where the selection of strategies is determined by a learning model. Prior theoretical analysis and experimental data show that after large number of plays the average number of agents who decide to enter, per round of the game, approaches the market capacity and, after a longer wait, agents are being sorted into two groups: the agents in one group rarely enter the market, and in the other, the agents enter almost all the time.  In this paper we obtain estimates of the characteristic times it takes for both patterns to emerge in the repeated plays of the game. The estimates are given in terms of the parameters of the game, assuming that the number of agents is large, the number of rounds of the game  per unit of time is large, and the characteristic change of the propensity per game is small. Our approach is based on the analysis of the partial differential equation for the function $f(t,q)$ that describes  the distribution of agents according to their level of propensity to enter the market, $q,$ at time $t.$

\end{abstract}


\end{frontmatter}

\section{Introduction}
\label{Intro}

A class of games in the study of social and econonmic behavior, called market entry games, describes a conflict  situation when players (agents) in a group choose between two strategies: enter the market or stay out, and the agent's payoff is determined solely by the number of agents who decide to enter and the action he or she takes.

The game has a single symmetric, mixed equilibrium and, a number  of asymmetric pure and mixed equilibria.

The extensive theoretical and experimental work has been done on understanding which, if any, of the equilibrium strategies emerge when the game is played repeatedly by not  agents who, independently from each other, try to adopt to changing ``market conditions.'' The situation can by formalized, by introducing into a model the individual propensities for agents to 
play a particular strategy. The propensities are updated after each round of the game. They might be  determined by the agent's payoffs, as in the {\it basic reinforcement learning model}, introduced in Erev and Roth (1998), or might depend on more information about the game available to agents, such as in the {\it fictitious stochastic play,} see  Fudenberg and Levine (1998).

When the number of players is large, the following patterns of behavior are typically observed and predicted by  learning models, see for example Duffy and Hopkins (2003):
\begin{enumerate}
\item {\it The average number of entries per round of the game quickly approaches the market capacity.} This is referred to as ``aggregate learning''.
\item {\it In a long-run of repeated plays agent's strategies converge to an asymmetric pure equilibrium, compatible with the market capacity.} This is called ``sorting''.
\end{enumerate}

Both phenomena are ubiquitous in the market entry games in which agents use either basic reinforcement learning or fictitious stochastic play. It is observed that the aggregate learning emerges quite quickly and it takes much longer time to observe sorting, see Duffy and Hopkins (2003).

The purpose of the present paper to give an estimate on the time scales of both phenomena, in terms of the number of agents, $N,$ the number of games, $M,$ played per unit of time, and the characteristic payoff $h$ per game. We show that the time of the aggregate learning is of the order
\[
\tau_{al}{}={}\frac{1}{MNh},
\]
and the time after which the sorting becomes noticeable is
\[
\tau_s{}={}\frac{1}{MNh^2}.
\]
One might expect that the formulas like these are only valid in a certain asymptotic regime, since the game is stochastic. This is indeed the case as the estimates are derived under the conditions that $N,M$ are large, $h$ is small, and $MNh$ is finite. Interestingly, the estimates are
the same for both models of the basic reinforcement learning and fictitious stochastic play.

Our approach is based on the derivation of a partial differential equation for the the distribution of agents among the propensity line. The equation is a drift-diffusion equation, with the drift velocity proportional to $(\tau_{al})^{-1},$ and the diffusion coefficient proportional to $(\tau_s)^{-1}.$

\section{The game and adaptive learning models}

There are $N$ agents participating in the game. Let $\delta^i$ denote the indicator function for agent $i$ to enter the game: $\delta^i=1$ if the agent enters, and $\delta^i=0,$ otherwise. Let $c\in \mathbb{N}$ be the capacity of the market, that we take for the simplicity of the presentation to be an integer. Let $m$ be the number of the agents who enter the market, $h$ be the characteristic payoff, and $v>0$ be the compensation for participating in the game. Then, the payoff to agent $i,$ can be defined, for example, as
\[
\pi^i{}={}\left\{
\begin{array}{ll}
v & \mbox{if } \delta^i=0,\\
v+h(c-m) & \mbox{if } \delta^i=1,
\end{array}
\right.
\]
see Erev and Rapoport (1998).

In the basic reinforcement learning, due to Erev and Roth (1998), the game is played repeatedly, and the state of the agent $i,$ is defined by the propensities to enter and stay out after $n^{th}$ round of the game: 
\[
(q^i_{1,n},q^i_{2,n})\in\mathbb{R}^2,
\]
The probability that agent $i$ uses in deciding to enter is given by
\[
y^i_n{}={}\frac{q^i_{1,n}}{q^i_{1,n}+q^i_{2,n}}.
\]
To reduce the number of parameters, in order to simplify the presentation, let us assume that $v=0.$ In this case 
the propensity to stay out does not changes in time: $q^i_{2,n}{}={}q^i_{2,0}.$ Furthermore, let us assume that for any $i=1..N,$ and some $q\in \mathbb{R},$ the propensity to stay our are the same for all agents:
\[
q^i_{2,n}{}={}q_0>0.
\]
Consiquently we need to consider only one the propensity to enter the market which we denote by
$
q^i_n.
$
Under such assumptions the probability for agent $i$ to enter the marker equals
\[
y^i_n{}={}\frac{q^i_n}{q^i_n+q}.
\]
To use this formula one has to make sure that propensities $q^i_n$ stay nonnegative. In fact, the explicit formula for the probability function will not be needed in our analysis, and we opt to use a generic probability function
\begin{equation}
\label{probability_function}
y^i_n{}={}p(q^i_{n}),
\end{equation}
where $p=p(q)\in[0,1]$ is strictly increasing, twice differentiable function such that
\[
p(-\infty){}{}={}1-p(+\infty){}={}0.
\]
In this way, the nonnegativty of the propensities is not required.

We consider two models of learning. In the model of basic reinforcement (with $v=0$), by Erev and Roth (1998), the propensity is increased/decreased by the amount of the payoff in $(n+1)^{th}$ game:
\begin{equation}
\label{SimpleLearning}
q^i_{n+1}{}={}q^i_n{}+{}h\delta^i_n(c-m_n),
\end{equation}
where $\delta^i_n$ is the indicator function of the action of player $i$ in $n^{th}$ game, and $m_n$ is the number of agent who enter the game.

In second model, the agents have more information about the game, which is relfected by fact that the propensity in \eqref{SimpleLearning} is increased by the
amount of the payoff agent $i$ would get if he/she played the opposite strategy:
\begin{equation}
\label{FictLearning}
q^i_{n+1}{}={}q^i_n+ h(c-m_n){}-{}h(1-\delta^i_n),
\end{equation}
see Duffy and Hopkins (2005).

\subsection{The method of the distribution function}
\label{df}

Using the theory of stochastic approximation of Bena\"{i}m (1999), Duffy and Hopkins (2005) prove that the repeated market entry games with either basic reinforcement learning or fictitious stochastic play, under rather generic condition, the agents strategies converge with probability one, to an asymmetric pure strategy equilibrium. In that approach, models \eqref{SimpleLearning} and \eqref{FictLearning} are considered as dynamical systems of size $N,$ that describe the individual propensity of all agents, as they involve under the stochastic updating rule. It is quite remarkable that the asymptotic behavior can be established for  such complicated systems.

In this paper we take a different approach, which is based on the derivation of the kinetic drift-diffusion equation for the distribution of the agents according to their propensity levels. 

Define the time step $\tau{}={}1/M,$ where $M$ is the number of rounds of the game per unit of time. The game takes place at times
\[
t_n{}={}n\tau, \quad n=1,2,3...
\]
If the initial propensities $q^i_0$ are discretized to the mesh $\{q_k{}={}kh\},$ $k\in\mathbb{Z},$ then for all times $t_n,$ propensities $q^i_n$ belong to the same mesh.

We are interested in the function $f(t_n,q)$
which is determined as the proportion of all agents that have propensity $q=q_k,$ $k\in\mathbb{Z},$ at time $t_n.$ That is, $f(t_n,q)$ is PMF (probability mass function) for the propensity of a randomly selected agent. We may write
\[
f(t_n,q){}={}\sum_k \alpha^n_k\delta(q-q_k),
\]
where $\delta(q-q_k)$ is the delta mass supported at $q_k,$
and $\alpha^n_k$ are non-negative numbers, summing up over $k$ to 1. They are defined as
\[
\alpha^n_k{}={}\frac{\mbox{\# of agents at time $t_n$ with propensity $q_k$}}{N}.
\]
We are interested in two integrals of $f.$ The first,
\begin{equation}
\label{a}
a(t_n){}={}\int p(q)f(t_n,q)\,dq
\end{equation}
is the fraction of the average number of entries to the market at $(n+1)^{th}$ round, and the second
\begin{equation}
\label{b}
 b(t_n){}={}\int p(q)(1-p(q))f(t_n,q)\,dq\geq 0,
\end{equation}
that we call the {\it coefficient of sorting}. 
The sorting of the population of agents into two groups is expressed by the smallness of $b(t),$ since it implies that  $f(t,q)$ is supported either on large negative $q's$ (rarely enter the market) or on large positive values (enter almost all the time).

We investigate the conditions under which $a(t)$ approaches the fraction of the market capacity $c/N$ and $b(t)$ converges to zero.  Thus, when working with the distribution function, we can not say to which particular equilibrium the system converges, but we still have enough information to say that the system does approaches an equilibrium and the equilibrium is a pure asymmetric one.

Let us also mention that studying distribution functions, instead of the dynamics of individual particles (agents) is a classical approach in Science, with the examples ranging from  the Boltzmann equation of gas dynamics and the diffusion processes  describing the Brownian motion to equations for distribution of commodities in social and economic studies, see Feller (1957), Ch. XIV, and  Pareschi \& Toscani (2014).


\begin{subsection}{Time scales}
We will show in \ref{Kinetic_equation} that in the asymptotic regime
\[
N\to\infty,\quad h\to0,\quad \tau\to0,\quad \frac{Nh}{\tau}\to r,
\]
for some $r\in\mathbb{R}^+,$ the density $f(t,q)$ of the basic reinforcement learning verifies the following nonlinear drift-diffusion equation. 
\begin{equation}
\label{eq:A}
\partial_t f+r(\kappa-a)\partial_q(p(q)f){}-{}\frac{1}{2}(   rNh(\kappa-a)^2{}+{}rh b)\partial_q^2 (p(q)f){}={}0,
\end{equation}
with
\[
a{}={}\int p(q)f(t,q)\,dq,\quad 
b{}={}\int p(1-p)f(t,q)\,dq,
\]
where $\kappa=c/N>0$ is the capacity of the market expressed as a fraction of the total population of agents.


The equation \eqref{eq:A} has rather simple structure.

It consists of diffusion of the density, with the diffusion coefficient
\[
\mu(t,q){}={}\frac{1}{2}( rNh(\kappa-a)^2+rhb)p'(q)>0,
\]
and the drift (transport) of the density $f$ with the velocity 
\[
v(t,q){}={}r(\kappa-a(t))p(q){}-{}\mu p'(q).
\]
In this notation equation \eqref{eq:A} takes the form:
\[
\partial_tf {}+{}\partial_q(v f){}-{}\partial_q(\mu\partial_q f){}={}0.
\]


We will show in \ref{Time_scales} that in the long run $a(t)$ approaches $\kappa$ and the characteristic time of this convergence is 
\begin{equation}
\label{T1}
\tau_{al}{}={}\frac{1}{r}.
\end{equation}

On the other hand $b(t)$ approaches zero and with the characteristic time
\begin{equation}
\label{T2}
\tau_s{}={}\frac{1}{rh}.
\end{equation}
Note in  particular that the sorting is slow compared
with the aggregate learning:
\[
\tau_s\gg \tau_{al}.
\]

The equation \eqref{eq:A} provides a convenient description  of processes governing the dynamics of  repeated games. Starting from  the initial distribution of propensities, the system quickly moves towards the state of the aggregate learning by punishing or rewarding {\it all} agents for deviations from the market capacity $\kappa.$ This is expressed in \eqref{eq:A} by the drift velocity being proportional to $(\kappa -a(t)).$

The stochastic nature of the decisions that agents make and their independence result in the tendency of the individual propensities to
spread out over the propensity space, which is expressed by the diffusion part of the equation \eqref{eq:A}. This results in the convergence of $f(t,q){}\to{}0,$ for every $q,$ meaning that fewer agents are using mixed strategies to play the game, leading to sorting. The sorting requires significant amount of time since diffusion coefficient is small.

\end{subsection}

\begin{subsection}{Time scales for the second model of learning}

For the model of learning \eqref{FictLearning} the density $f$ verifies a similar equation
\begin{eqnarray}
\label{eq:B}
\partial_t f{}+{}r(\kappa-a)\partial_qf{}-{}\frac{1}{2}(   rNh(\kappa-a)^2{}+{}rh b)\partial_q^2 f{}={}0,
\end{eqnarray}
which lead to the same estimates for time scales \eqref{T1} and \eqref{T2}.

Note here the difference between equations \eqref{eq:A} and \eqref{eq:B}: the drift velocity and the diffusion coefficient in the model of fictitious play is homogeneous (independent of) in the propensity $q,$ whereas in \eqref{eq:A} both proportional to the probability $p(q).$ The later case indicates that the intensity of the drift and diffusion is larger for agents who have higher probability (and thus propensity) to enter. This, in fact, is expected, as the only way to change the propensity of an agent in the reinforcement model is for her to enter the market.  

This difference results in the sorting being somewhat faster for agents with low propensity in the model  \eqref{eq:B}, once the system gets close to the equilibrium. 
\end{subsection}

\appendix

\section{The drift-diffusion equation}
\label{Kinetic_equation}
In this section we give a heuristic derivation of drift-diffusion equations \eqref{eq:A} and \eqref{eq:B}. Consider a repeated market entry game with basic reinforcement learning \eqref{probability_function}, \eqref{SimpleLearning}.

Let $f(t_n,q)$ be the PMF of the propensity of a randomly selected agent, as in section \ref{df}, and $X_n$ be a random variable with this PMF.

At time $t=t_n$ we select an agent at random and observe her propensity level: $\hat{q}{}={}X_n.$ The probability that the agent will enter market at the next round is $\hat{p}{}={}p(\hat{q}).$ $\hatq$ belongs to the mesh $\{kh\}$ and we denote the corresponding mesh number
\[
\hat{k}{}={}\frac{\hat{q}}{h}\in\mathbb{Z}.
\]

If $k\not=\hatk,$ there are $Nf(t_n,q_k){}={}N\alpha^n_k$ agents that have propensity $q_k;$ if $k=\hatk,$
there are $N\alpha^n_{\hatk}-1$ agents (besides the one we selected) that have propensity $\hatq.$

Since each agent behaves independently from others, the number of agents among $N\alpha^n_k$ (or $N\alpha^n_{\hatk}-1$) who will enter at the next round is a binomial random variable,
that we denote $X^n_k$ (or $X^n_{\hatk}$):
\[
X^n_k\in B(p_k,N\alpha^n_k),\quad p_k=p(q_k),\,k\not=\hatk,
\]
and 
\[
X^n_{\hatk}\in B(p_{\hatk},N\alpha^n_{\hatk}-1).
\]
Here $B(n,p)$ stands for the  binomial distribution of successes in $n$ Bernoulli trials, with the probability of success $p.$

We compute the expectation and the variance of $X^n_k$,
\[
E[X^n_k]{}={}N\alpha^n_kp_k,\,V(X^n_k){}={}N\alpha^n_kp_k(1-p_k),\quad k\not=\hatk,
\]
and
\[ 
E[X^n_{\hatk}]{}={}(N\alpha^n_{\hatk}-1)p_{\hatk},\,V(X^n_{\hatk}){}={}(N\alpha^n_{\hatk}-1)p_{\hatk}(1-p_{\hatk}).
\]
Let $\hatm$ be the total number of agents who enter the market, not counting the selected one, i.e.,
\[
\hatm{}={}\sum_k X^n_k.
\]
Then, using the pairwise independence of $X^n_k$'s,
\[
E[\hatm]{}={}N\int p(q)f(t_n,q)\,dq{}-{}p_{\hatk},
\]
and
\[
V(\hatm){}={}N\int p(q)(1-p(q))f(t_n,q)\,dq{}-{}p_{\hatk}(1-p_{\hatk}).
\]
Let $\hatdelta$ be the indicator function for the selected agent to enter the market, which is,  by design of the model, independent of $\hatm.$ For such random variable,
\[
\mbox{Prob}(\hatdelta=1){}={}1-\mbox{Prob}(\hatdelta=0){}={}p_{\hatk}.
\]
At the next time, $t=t_{n+1},$ the propensity of the selected player to enter the market by \eqref{SimpleLearning} equals
\begin{equation}
\label{Xn}
X_{n+1}{}={}X_n+\hatdelta h( c- \hatm -1).
\end{equation}
{\it  We will assume that $X_{n+1}$ is a good approximation of the propensity to enter the market of a randomly selected agent at time $t=t_{n+1}.$ In other words, the distribution of $X_{n+1}$ is given by $f(t_{n+1},q).$}

This assumption suffices to the derive the equation for $f(t,q).$ Let $\phi$ be a test function and compute
\begin{multline}
\int \phi(\hatq)f(t_{n+1},\hatq)\,d\hatq{}={}E[\phi(X_{n+1})]
{}={} E[\phi(X_n+\hatdelta h( c- \hatm )]\\
{}={}
\int \left(p(\hatq)E[\phi(\hatq+h(c-\hatm -1))]+(1-p(\hatq))\phi(\hatq)  \right)f(t_n,\hatq)\,d\hatq.
\end{multline}

Using the Taylor's expansion
\[
\phi(\hatq+z){}={}\phi(\hatq){}+{}\phi'(\hatq)z{}+{}\frac{\phi''(\hatq)}{2}z^2 + o(z^2),
\]
we obtain from the previous computation:
\begin{multline}
\label{K2}
\int \phi(\hatq)\left( \frac{f(t_{n+1},\hatq)-f(t_n,\hatq)}{\tau}\right)\,d\hatq{}={}
\frac{h}{\tau}\int \phi'(\hatq)p(\hatq)E[c-\hatm-1]f(t_n,\hatq)\,d\hatq\\
{}+{}
\frac{h^2}{2\tau}\int\phi''(\hatq) p(\hatq)E[(c-\hatm-1)^2]f(t_n,\hatq)\,dq{}+{}\frac{ho(h)}{\tau}. 
\end{multline}
Next we compute
\begin{multline}
E[c-\hatm-1]{}={}c- N\int p(q)f(t_n,q)\,dq{}+{}p(\hatq)-1\\
\approx N\left(c/N-\int p f(t_n,q)\,dq\right),
\end{multline}
where we assumed that $N$ is large.
In a similar way,
\begin{multline}
E[(c-\hatm-1)^2]{}={}
(c-E[\hatm]-1)^2{}+{} V(\hatm)\\{}\approx{}
\left( N^2(c/N -\int p f(t_n,q)\,dq)^2{}+{}N(\int p(1-p)f(t_n,q)\,dq)\right).
\end{multline}

Finally we assume that $
\tau,\, Nh$ are small and $r=\dfrac{Nh}{\tau},\,\kappa=\dfrac{c}{N} $ are finite. 

Returning to \eqref{K2} and retaining only higher order terms we obtain an integral equation
\begin{multline}
\int \phi(q)\partial_t f(t,q)\,dq{}={}
r\int \phi'(q)p(q)(\kappa -\int p f\,dq)f\,dq\\
{}+{}
\frac{rNh}{2}\int \phi''(q) p(q)\left( (\kappa -\int p f\,dq)^2{}+{}\int p(1-p)f\,dq/N   \right)f\,dq. 
\end{multline}
Since the equation holds for an arbitrary test function $\phi,$ we obtain a partial differential equation for $f:$
\[
\partial_t f{}+{}r(\kappa-a)\partial_q(pf){}-{} \frac{1}{2}(   rNh(\kappa-a)^2{}+{}rh b)\partial_q^2(pf){}={}0,
\]
where
\[
a{}={}\int p(q)f(t,q)\,dq,\, b{}={}\int p(q)(1-p(q))f(t,q)\,dq,\, \kappa{}={}\frac{c}{N},\,r{}={}\frac{Nh}{\tau}.
\]

For the second model of learning \eqref{FictLearning}, instead of \eqref{Xn}, the updating rule prescribes
\[
X_{n+1}{}={}X_n+ h(c-\hatm){}-{}\hatdelta h.
\]

Repeating the arguments of the previous case  we obtain the equation
\[
\partial_t f{}+{}r(\kappa-a)\partial_q(f){}-{}\frac{1}{2}(   rNh(\kappa-a)^2{}+{}rh b)\partial_q^2(f){}={}0,
\]
where $a,b$ are the same as above.

\section{Two time scales}
\label{Time_scales}
Formulas \eqref{T1} and \eqref{T2} are obtained by taking moments of the equation \eqref{eq:A}. We multiply the equation by $p(q)$ and integrate in $q.$ Assuming the $f(t,q)$ decays fast enough at infinity we find an ODE for the average entry rate $a(t):$
\begin{equation}
\label{eq:a1}
\frac{da}{dt}{}={}r\left(\int p'pf(t,q)\,dq\right)(\kappa -a)
{}+{}\frac{1}{2}(rNh(\kappa-a)^2{}+{}rh b)\int p''pf\,dq.
\end{equation}
The ODE is still not closed as it depends on an integral of $f.$ However, since $p(q),\,p'(q)>0$ and 
\[
\int p' pf\,dq{}\leq{}\max\{p'(q)p(q)\},
\]
 it is of the order 1 and we substitute it with a suitable positive constant $c(p).$  Moreover, since $Nh\approx \tau$ and $h$ are small, the contribution of the last term in \eqref{eq:a1} can be ignored and we obtain the equation
\begin{equation}
\label{eq:a}
\frac{da}{dt}{}={}rc(p)\left(\kappa -a\right),
\end{equation}
and 
\[
a(t){}={}\kappa + (a(0)-\kappa)e^{-c(p)rt}.
\]
Thus, the ratio $(a(t)-\kappa)/(a(0)-\kappa)$ decreases to zero with the characteristic time length $\tau_{al}{}={}1/r.$

It is harder to obtain an analytical expression for the coefficient of sorting. Qualitatively, the sorting occurs due to the diffusion of the density $f.$ We will proceed heuristically, postulating that the rate of decrease of the coefficient of sorting is proportional to the diffusion
\[
(rNh(\kappa-a(t))^2{}+{}rh b(t))/2.
\]
Since this quantity is asymptotically smaller the the drift $r(\kappa-a(t)),$  the diffusion is of the order
$
rhb(t)/2.
$
Thus, we obtain that
\[
\frac{db}{dt}\approx -\frac{rh}{2}b,
\]
and $b(t)\approx  b(0)e^{-rht/2}.$ This formula implies that the ratio $b(t)/b(0)$ decreases to zero at the characteristic time $\tau_s{}={}2/(rh).$


\begin{thebibliography}{00}


\bibitem{Benaim} Bena\"{i}m,~M., 1999. {\sl Dynamics of stochastic algorithms.\/} S\'{e}minaire de probabilit\'{e}s 
Strasbourg, 33, 1--68.

\bibitem{DuffyHopkins} Duffy,~J., Hopkins,~E., 2005. {\sl Learning, information, and sorting in market entry games: theory and evidence.\/} Games  Econ. Behav. 51, 31--62. 

\bibitem{ErevRapoport} Erev,~I., Rapoport,~A., 1998. {\sl Coordination, ``magic'', and reinforcement learning in a market entry game.\/} Games Econ. Behav. 23, 46--175.

\bibitem{ErevRoth} Erev,~I, Roth,~A.~E., 1998. {\sl Predicting how people play games: reinforcement learning in experimental games with unique, mixed strategy equilibrium.\/} American Econ. Review 88, 848--881.

\bibitem{Feller} Feller,~W., 1957. {\sl An Introduction to Probability Theory and Its Applications.\/} John Wiley \& Sons, New York, NY.

\bibitem{FudenbergLevine} Fudenberg,~D., Levine,~D.K.,  1998. {\sl The Theory of Learning in Games.\/} MIT Press, Cambridge, MA.


\bibitem{ParschiToscani} Pareschi,~L., Toscani,~G., 2014. {\sl Interacting Multiagent Systems.\/} Oxford University Press, Oxford, UK.


\end{thebibliography}
\end{document}